\documentclass[a4paper,11pt]{amsart}
\usepackage{amsfonts}
 \usepackage{enumerate}
 \usepackage{amsmath, amsthm, amscd, amssymb}
\usepackage{pdfsync}
\usepackage[title]{appendix}
\usepackage[all]{xy}
\usepackage{latexsym,graphicx}
\usepackage[latin1]{inputenc}
\usepackage{xspace}
\usepackage{array}
\usepackage{newclude}
\usepackage{hyperref}
\usepackage{url}
\usepackage{enumitem}
\hypersetup{
  colorlinks,
  citecolor=black,
  linkcolor=magenta,
  urlcolor=black}
  
\renewcommand*{\HyperDestNameFilter}[1]{\jobname-#1} 
\numberwithin{equation}{section}
\usepackage[centering, marginparwidth=2cm]{geometry}
\usepackage{marginnote}

\addtocontents{toc}{\setcounter{tocdepth}{1}}

\usepackage[nameinlink]{cleveref}

\usepackage{tikz-cd}
\usetikzlibrary{cd}
\usepackage{leftidx}
\usetikzlibrary{positioning}
\usepackage{dsfont}
\usepackage{tikz}
\usetikzlibrary{matrix,arrows,decorations.pathmorphing,positioning}

\newcommand\blfootnote[1]{%
\begingroup
\renewcommand\thefootnote{}\footnote{#1}%
\addtocounter{footnote}{-1}%
\endgroup
}
  
  \newcommand{\Addresses}{{
  \bigskip
  \footnotesize

\textsc{I.H.E.S., Universit\'e Paris-Saclay, CNRS, Laboratoire Alexandre
  Grothendieck. 35 Route de Chartres, 91440 Bures-sur-Yvette
  (France)}\par\nopagebreak
  \textit{E-mail address}, G.~Baldi: \texttt{baldi@ihes.fr}  

  \medskip

  \textsc{Humboldt Universit\"{a}t zu Berlin (Germany)}
  \par\nopagebreak
  \textit{E-mail address}, B.~Klingler: \texttt{bruno.klingler@hu-berlin.de}
  
  \medskip
  
\textsc{I.H.E.S., Universit\'e Paris-Saclay, CNRS, Laboratoire Alexandre Grothendieck. 35 Route de Chartres, 91440 Bures-sur-Yvette (France)}\par\nopagebreak
  \textit{E-mail address}, E.~Ullmo: \texttt{ullmo@ihes.fr}
}}

\theoremstyle{plain}
\newtheorem{theor}{Theorem}[section]

\theoremstyle{definition}

\newtheorem{rmk}[theor]{Remark}

\theoremstyle{remark}
\numberwithin{equation}{subsection}

\newcommand{\CC}{\mathbb{C}}
\newcommand{\NN}{\mathbb{N}}
\newcommand{\ZZ}{\mathbb{Z}}
\newcommand{\QQ}{\mathbb{Q}}

\newcommand{\an}{\textnormal{an}}

\newcommand{\ad}{\textnormal{ad}}

\newcommand{\Zar}{\textnormal{Zar}}

\newcommand{\PP}{\mathbf{P}}

\newcommand{\VV}{{\mathbb V}}

\newcommand{\HL}{\textnormal{HL}}
\newcommand{\atyp}{\textnormal{atyp}}

\newcommand{\pos}{\textnormal{pos}}

\newcommand{\N}{\mathbb{N}}

\newcommand{\prim}{\textnormal{prim}}

\newcommand{\Z}{\mathbb{Z}}
\newcommand{\Q}{\mathbb{Q}}

\newcommand{\Oo}{\mathcal{O}}
\newcommand{\oQ}{\overline{\QQ}}

\usepackage[normalem]{ulem}

\newcommand{\C}{\mathbb{C}}

\newcommand{\Qbar}{\overline{\mathbb{Q}}}

\newcommand{\PGL}{\mathbf{PGL}}
\newcommand{\G}{{\mathbf G}}

\newcommand{\cM}{\mathcal{M}}

\newcommand{\cO}{\mathcal{O}}

\usepackage{microtype}
\usepackage{fullpage}

\begin{document}

\title{On the Geometric Zilber--Pink Theorem and the Lawrence--Venkatesh method}\date{\today}  
\author{Gregorio Baldi, Bruno Klingler, and Emmanuel Ullmo} 

\begin{abstract} Using our recent results on the algebraicity of the
  Hodge locus for variations of Hodge structures of level at least
  $3$, we improve the results of Lawrence-Venkatesh in direction of the
  refined Bombieri--Lang conjecture. 
\end{abstract}
\maketitle

\section{Results}
The aim of this short note is to explain how the \emph{Geometric
  Zilber--Pink conjecture} for variations of Hodge structures, recently established by the authors
\cite{2021arXiv210708838B}, can be used to improve the main result of
Lawrence and Venkatesh \cite{zbMATH07233321}, giving a special case of
the \emph{refined Bombieri--Lang conjecture}.

\medskip
Let $U_{n,d}$ be the Hilbert scheme of smooth hypersurfaces of degree $d$ in
$\PP^{n+1}$, this is a smooth scheme over $\ZZ$. In a recent
breakthrough, Lawrence and Venkatesh proved the following: 

\begin{theor}[{\cite[Theorem 10.1, Proposition 10.2]{zbMATH07233321}}]\label{lv}
There exist $n_0 \in \NN_{\geq 3}$ and a function $d_0: \NN \to \NN$ such that,
\begin{equation}\label{conditionLV}
\textnormal{ for every} \; n \geq n_0 \text{   and   } d\geq d_0(n),
\end{equation}
the set $U_{n,d}(\Z[S^{-1}])$ is not Zariski dense in 
$U_{n,d}$, for every finite set of primes $S$.
\end{theor}

Consider $f_{n,d}: X_{n,d} \to U_{n,d}$
the universal family of smooth degree $d$ hypersurfaces in
$\PP^{n+1}$. We denote by $\VV$ the polarized $\ZZ$-variation of
Hodge structure $(R^n {f ^\an_{n, d, \CC}}_* \ZZ)_\prim$ on $U_{n, d,
  \CC}$ and by $\Phi: U_{n,d, \CC}^\an \to \Gamma \backslash D$ the
associated period map. An irreducible algebraic subvariety $Y \subset
S$ is said to be {\em of positive period dimension} if $\Phi(Y_\CC^\an)$
has positive dimension. We prove the following reinforcement of \Cref{lv}:

\begin{theor}\label{maincor} 
As long as \eqref{conditionLV} is satisfied, there exists a closed strict subscheme 
$E\subset U_{n,d}$ such that, for all finite set of primes $S$, we have
\begin{displaymath}
\overline{U_{n,d}(\Z[S^{-1}])}_{\pos} \subset E,
\end{displaymath}
where $\overline{U_{n,d}(\Z[S^{-1}])}_{\pos}$ denotes the union of the irreducible
components of the Zariski closure of $U_{n,d}(\Z[S^{-1}])$ in
$U_{n,d}$ of positive period dimension.
That is: the Zariski closure of $U_{n,d}(\Z[S^{-1}]) - E(\Z[S^{-1}])$ has period dimension zero.
\end{theor}

\begin{rmk}\label{lastrmk}
  The complement of $U_{n,d}$ in $\PP^{N(n, d)}$ is a
  hypersurface. Hence $U_{n,d}$ is an open affine subvariety, stable
  under the natural $\PGL(n+2)$-action on $\PP^{N(n, d)}$. Let
  $\cM_{n, d}:= [\PGL(n+2) \backslash U_{n,d}]$ be the stack of smooth
  hypersurfaces in $\PP^{n+1}$ of degree $d$. This is a finite type
  separated Deligne-Mumford algebraic stack over $\ZZ$ with affine
  coarse space, see \cite{MR3022710}.  The period map
  $\Phi: U_{n, d, \CC}^\an  \to \Gamma \backslash D$ factorizes through
  $\cM_{n, d, \CC} ^\an$. Notice moreover
  that the Torelli theorem (see e.g. \cite{2020arXiv200409310V} and
  references therein) assures 
    that the period map $\cM_{n, d, \CC} ^\an \to \Gamma \backslash D$
    is quasi-finite for $(n, d) \not= (2,3)$.

The moduli stack $\cM_{n, d, \CC}^\an$ is thus
\emph{Brody hyperbolic}. A famous conjecture of
Bombieri-Lang (see for instance \cite[Chapter F.5.2]{zbMATH01466163})
thus predicts that $\cM_{n,
  d}(\ZZ[S^{-1}])$ is finite. In particular $E$ in
\Cref{maincor} should be empty. 
\end{rmk}




We first recall in \Cref{Hodge locus} (a special case of) the
Geometric Zilber--Pink conjecture mentioned above, which is a purely
geometric result; and then in \Cref{integralpoints} the Lawrence-Venkatesh
method, which is of arithmetic nature. In \Cref{Main} we explain what can be obtained by combining the two results.

\subsection*{Acknowledgements}We thank A. Javanpeykar for
  pointing out some inaccuracies in a previous draft, B. Lawrence and
  A. Venkatesh for their interest. B.K. is partially supported by the
  Grant ERC-2020-ADG n.101020009 - TameHodge.

\section{The geometry of the Hodge locus} \label{Hodge locus}
Let $f: X \to S$ be a smooth projective morphism of smooth irreducible complex
quasi-projective varieties, of relative dimension $n$. The primitive Betti
cohomology $H^n(X_s^\an, \ZZ)_\prim$ of the fibres $X_s$, $s \in S(\CC)$,
form a polarized $\ZZ$-variation of Hodge structures $\VV$ on the complex manifold
$S^\an$, described by a complex analytic period map $\Phi: S^\an \to \Gamma \backslash
D$ (we refer for instance to \cite{2021arXiv210708838B} for more details on period maps).  Motivated by the study of the Hodge conjecture for the
fibres of $f$, one defines the Hodge locus $\HL(S, \VV^\otimes)$ as the locus of
points $s\in S^{\an}$ for which the Hodge structure
$H^n(X_s^{\an},\Q)_\prim$ admits more \emph{Hodge
tensors} than the primitive cohomology of the very
general fibre. Here a Hodge class of a pure $\ZZ$-Hodge structure $V=(V_\ZZ, F^\bullet)$ is a class in $V_\Q$
whose image in $V_\CC$ lies in the zeroth piece $F^0 V_\CC$ of the
Hodge filtration, or equivalently a morphism of Hodge 
structures $\Q(0) \to V_\QQ$; and a Hodge tensor for $V$ is a Hodge class
in $V^\otimes:= \bigoplus_{a,b\in
  \N} V^{\otimes a} \otimes( V^\vee)^{\otimes b}$, where
$V^\vee$ denotes the Hodge structure dual to $V$. Cattani, Deligne and
Kaplan \cite[Theorem 1.1]{CDK95} proved in particular that the Hodge
locus $\HL(S, \VV^\otimes)$ is a \emph{countable} union of irreducible 
algebraic subvarieties of $S$, called the special subvarieties of $S$
for $\VV$ (or $f$). We denote by $\HL(S, \VV^\otimes)_{\pos}$ the Hodge locus \emph{of positive period dimension}, that is the 
union of the special subvarieties whose image under $\Phi$ has
positive dimension in $\Gamma \backslash D$.

Let $\PP^{N(n, d)}_\QQ:= \PP (H^0(\PP^{n+1}_\QQ,
\cO_{\PP^{n+1}_\QQ}(d)))$ be the Hilbert scheme of hypersurfaces $X$ of
$\PP^{n+1}_\QQ$ of degree $d$ (where $N(n,
d)=\binom{n+d+1}{d}-1$). Let $U_{n, d} \subset \PP^{N(n, d)}_\QQ$ be 
the Zariski-open Hilbert scheme of smooth 
hypersurfaces $X$ and consider 
$
f_{n,d}: X_{n,d} \to U_{n,d}
$
the universal family of smooth degree $d$ hypersurfaces in
$\PP^{n+1}_\QQ$. We denote by $\VV$ the polarized $\ZZ$-variation of
Hodge structure $(R^n {f ^\an_{n, d, \CC}}_* \ZZ)_\prim$ on $U_{n, d,
  \CC}$. We write $\HL(U_{n,d, \CC}, \VV^\otimes)$ for its Hodge locus and
$\HL(U_{n,d, \CC}, \VV^\otimes)_{\pos}$ for its Hodge locus of positive period
dimension.

In our previous paper we have established the following 
as a particular case of our main result:

\begin{theor}[{\cite[Corollary
    2.7]{2021arXiv210708838B}}] \label{hypersurface}

If $n\geq 3, d\geq 5$ and $(n,d)\neq (4,5)$ then the
Hodge locus $\HL(U_{n, d, \CC}, \VV)_\pos$ of positive period dimension is
a closed (not necessarily irreducible) algebraic subvariety of $U_{n,d, \CC}$. That is, there are only finitely many (rather than countably many)
maximal strict special subvarieties of $U_{n, d, \CC}$ for $\VV$ of
positive period dimension.
\end{theor}

\begin{rmk}\label{rmk11}
For what follows, the
easier \cite[Theorem 5.1]{2021arXiv210708838B} would actually be
enough (that is the \emph{Geometric Part of Zilber-Pink}, for
weakly-special subvarieties). 
\end{rmk}

\section{Non-denseness of integral points and the Lawrence--Venkatesh
  method}\label{integralpoints}



The following elucidation of the Lawrence-Venkatesh method for proving
\Cref{lv} will be crucial for us. Lawrence and Venkatesh actually
prove that (quoting the third
paragraph of \cite[Section 1.1]{zbMATH07233321}) \emph{the monodromy
  for the universal family of hypersurfaces must drop over each
  component of the Zariski closure of the integral points} (see also
the last three lines of \cite[Theorem 10.1]{zbMATH07233321}):  for any $S$, there
exists a closed subscheme $V_{S}$ of
$U_{n,d}$ (over $\Z[S^{-1}]$) whose irreducible components are of
  positive period dimension and not \emph{monodromy
  generic}, such that $U_{n,d}(\Z[S^{-1}])-V_{S}(\Z[S^{-1}])$ is
finite. By the Deligne-Andr\'{e} monodromy theorem (see
for example \cite[Section 3 and 4]{2021arXiv210708838B}) and the fact
that the $\Z$VHS $\VV$ is irreducible, it follows that each $V_{S}$
lies in the Hodge locus of positive
period dimension $\HL(U_{n,d},
\VV^\otimes)_{\pos}$.
\footnote{In fact, and to justify \Cref{rmk11}, their proof actually shows
  that each $V_{S}$ is contained in the \emph{atypical Hodge
    locus} of positive period dimension. Such subspace of the Hodge
  locus is proven to be non-Zariski dense in $U_{n,d}$ in
  \cite{2021arXiv210708838B} as a first step towards
  \Cref{hypersurface}, but it holds true for any variety supporting
  any variation of Hodge structures whose adjoint generic
    Mumford-Tate group is simple.}

\begin{rmk} \label{finalremark}
The Lawrence-Venkatesh method requires the choice
of an auxiliary prime number $p$, and the choice of an identification between $\C$
and $\overline{\Q}_p$. Indeed, to prove that the $\ZZ[S^{-1}]$-points of
$U_{n,d}$ are not Zariski dense, Lawrence and Venkatesh prove that some \emph{p-adic
  period map} sending $x \in U_{n,d}(\ZZ[S^{-1}])$ to some $p$-adic
representation of the absolute Galois group of $\Q_p$ has fibers that are
not Zariski dense in $U_{n,d}$. This is done by working on a residue
disk in $U_{n,d}(\Q_p)$ and the $p$-adic and complex period maps are
then related by a study of the Gauss-Manin connection \cite[Lemma
3.2]{zbMATH07233321}. What their proof actually shows, with respect
to our fixed embedding $\oQ \subset \CC$ is that for each $S$, there
exists an automorphism $\iota_p$ of $\CC$ such that $V_{S}$ is
contained in $\HL(U_{n, d}, \VV^\otimes)_\otimes^{\iota_p} \subset U_{n, d,
  \CC}$. What allows us to say that $V_{S}$ lies in $\HL(U_{n, d}, \VV^\otimes)_\pos$ is the fact that $\HL(U_{n,d}, \VV^\otimes)_\pos$ is actually
defined over some number field, as one sees by combining \Cref{hypersurface} and
\cite[Theorem 1.10]{2020arXiv201003359K}. 
\end{rmk}

\begin{rmk}\label{rmk1}
Let us emphasize that both  \cite[Theorem 10.1]{zbMATH07233321} and \Cref{hypersurface}
build on the Ax-Schanuel theorem 
\cite{MR3958791}, a deep and general theorem establishing strong
functional transcendence properties of period maps. Actually, in
Lawrence-Venkatesh, a $p$-adic version of such a result is used, see
indeed \cite[Lemma 9.3]{zbMATH07233321}. 
\end{rmk} 
 
 \section{Proof of the Main result} \label{Main}

The proof is essentially a combination of \Cref{lv}
and
\Cref{hypersurface}.

It follows from \Cref{hypersurface} that $\HL(U_{n,d}, \VV^\otimes)_\pos$ is a (closed, strict) algebraic
subvariety of $U_{n,d}$ and, thanks to the elucidation of \Cref{lv}, we have 
\begin{displaymath}
\bigcup_{S} \overline{U_{n,d}(\Z[S^{-1}]})_{\pos} = \bigcup_{S}
V_{S} \subset  \HL(U_{n,d}, \VV^\otimes)_{\pos}, 
\end{displaymath}
where the union ranges over
all finite set of primes $S$. It follows from \Cref{hypersurface} that 
\begin{displaymath}
E':= \overline{\bigcup_{S} V_{S} }^{\Zar} \subset \HL(U_{n,d}, \VV^\otimes)_\pos.
\end{displaymath}
We remark here that the above inclusion may happen to be
strict. Therefore we obtained a closed $\Qbar$-subvariety $E' \subset
\HL(U_{n,d}, \VV^\otimes)_{\pos}$ containing all $V_{S}$ (seen as
$\Qbar$-varieties). The Zariski closure $E$ in $\PP^N_\Z$ of $E'$
enjoys the desired property. The proof of
the Theorem is concluded.

\bibliographystyle{abbrv}

\bibliography{biblio.bib}

\Addresses

\end{document}